\magnification= \magstep1
\input amstex
\documentstyle{amsppt}
\define\defeq{\overset{\text{def}}\to=}

\define\Gal{\operatorname{Gal}}

\def \isom {\overset \sim \to \rightarrow}

\define\Spec{\operatorname{Spec}}
\def \et{\operatorname {et}}
\def \Ker{\operatorname {Ker}}
\def \id{\operatorname {id}}
\def \Fr{\operatorname {Fr}}
\def \char{\operatorname {char}}
\def \Spec{\operatorname {Spec}}
\def \Spf{\operatorname {Spf}}
\def \Sp{\operatorname {Sp}}
\def \geo{\operatorname {geo}}

\def \alg{\operatorname {alg}}
\def \cd{\operatorname {cd}}

\def \Ind{\operatorname {Ind}}
\def \rig{\operatorname {rig}}

\def \Aut{\operatorname {Aut}}
\def \Im{\operatorname {Im}}
\NoRunningHeads
\NoBlackBoxes
\topmatter

\title 
\'Etale fundamental groups of affinoid $p$-adic curves 
\endtitle
\bigskip

\author
Mohamed Sa\"\i di
\bigskip
\centerline 
{\it In memory of Si M'hamed, my father.}
\endauthor

\abstract We prove that the {\it geometric} \'etale fundamental group of a (geometrically connected) 
{\it rigid smooth $p$-adic affinoid curve} 
is a {\it semi-direct factor} of a certain profinite {\it free} group. We prove that the maximal pro-$p$ (resp. maximal prime-to-$p$) quotient
of this geometric \'etale fundamental group is pro-$p$ {\it free of infinite rank} (resp. (pro-)prime-to-$p$ {\it free of finite computable rank}). 
\endabstract

\toc

\subhead
\S0. Introduction/Main Results
\endsubhead

\subhead
\S1. Background
\endsubhead

\subhead
\S2. Geometric fundamental groups of annuli of thickness zero
\endsubhead

\subhead
\S3. Geometric fundamental groups of affinoid $p$-adic curves
\endsubhead

\endtoc

\endtopmatter

\document 

\subhead
\S 0. Introduction/Main Results
\endsubhead
A classical result in the theory of \'etale fundamental groups is the description of the structure of the geometric \'etale fundamental group of 
an affine, smooth, and geometrically connected curve over a field of characteristic $0$ (cf. [Grothendieck], Expos\'e XIII, Corollaire 2.12). 
In this paper we investigate the structure of the geometric \'etale fundamental group of a smooth affinoid $p$-adic curve.

Let $R$ be a complete discrete valuation ring, $K\defeq \Fr (R)$ the quotient field of $R$, and
$k$ its residue field which is algebraically closed of characteristic $p\ge 0$. 
Let $X_K$ be a smooth, proper, and geometrically connected rigid $K$-curve, $\Cal U\hookrightarrow X_K$ a $K$-{\it affinoid} rigid subspace
with $\Cal U$ geometrically connected and $X_K\setminus \Cal U$ is the disjoint union of $K$-rigid 
open unit discs $\{\Cal D_{i}^o\}_{i=1}^m$ with centres $\{x_i\}_{i=1}^m$, $x_i\in X_K(K)$ (cf. $\S3$ for more details, as well as Theorem 3.1 which asserts that 
any $K$-affinoid smooth curve can be embedded, after possibly a finite extension of $K$, into a proper and smooth rigid $K$-curve whose complement is as above).

Let $S\subset \Cal U$ be a (possibly empty) finite set of points and
$T\subset \bigcup_{i=1}^m \Cal D_{i}^o$ a finite set of points of 
$X_K$. (We also denote, when there is no risk of confusion, by $X_K$ the projective, smooth, and geometrically connected algebraic 
$K$-curve associated to the rigid curve $X_K$ via the rigid GAGA functor.) 
We have an exact sequence of \'etale fundamental groups
$$1\to \pi_1(X_{K}\setminus (T\cup S))^{\geo}\to \pi_1(X_K\setminus (T\cup S))\to 
\Gal (\overline K/K)\to 1,$$
where $\pi_1(X_K\setminus (T\cup S))$ is the arithmetic fundamental group of the (affine) curve $X_K\setminus (T\cup S)$,
and by passing to the projective limit over all finite sets of points 
$T\subset  \bigcup_{i=1}^m \Cal D_{i}^o$ we obtain an exact sequence
$$1\to \underset{T} \to {\varprojlim}\ \pi_1(X_{K}\setminus (T\cup S))^{\geo}
\to \underset {T} \to {\varprojlim}\ \pi_1(X_{K}\setminus (T\cup S))\to \Gal (\overline 
K/K)\to 1.$$
The profinite group $\underset{T} \to {\varprojlim}\ \pi_1(X_{K}\setminus (T\cup S))^{\geo}$
is {\it free} if $\char (K)=0$ as follows from the well-known structure of the geometric \'etale fundamental groups of (affine) curves
in characteristic zero (cf. loc. cit.).  Write $\pi_1(\Cal U\setminus S)^{\geo}$ for the 
geometric \'etale fundamental group (in the sense of Grothendieck) of $\Cal U\setminus S$ (cf. 2.1 for a precise Definition). 
One of our main results is the following (cf. Theorem 3.4, Proposition 3.5, and Theorem 3.7).

\proclaim {Theorem A}
Assume $\char (K)=0$ and $\char(k)=p>0$. 
Let $\ell$ be a prime integer (possibly equal to $p$). Then the morphism $\Cal U\to X_K$ induces (via the rigid GAGA functor) a continuous 
homomorphism $\pi_1(\Cal U\setminus S)^{\geo}
\to \underset{T} \to {\varprojlim}\  \pi_1(X_K\setminus (T\cup S))^{\geo}$ 
(resp. $\pi_1(\Cal U\setminus S)^{\geo,\ell}
\to \underset{T} \to {\varprojlim}\  \pi_1(X_K\setminus (T\cup S))^{\geo,\ell}$ between the maximal pro-$\ell$ quotients) 
which makes $\pi_1(\Cal U\setminus S)^{\geo}$ (resp. $\pi_1(\Cal U\setminus S)^{\geo,\ell}$) 
into a semi-direct factor of $\underset{T} \to {\varprojlim}\ \pi_1(X_K\setminus (T\cup S))^{\geo}$
(resp. a direct factor of $\underset{T} \to {\varprojlim}\ \pi_1(X_K\setminus (T\cup S))^{\geo,\ell}$)
(cf. Definitions 1.2 and 1.4 for the meaning of the terms direct factor and semi-direct factor).
Moreover, the pro-$\ell$ group $\pi_1(\Cal U\setminus S)^{\geo,\ell}$ is free of infinite rank if $\ell=p$, and of finite computable rank if $\ell\neq p$.
\endproclaim

In Theorem 3.3 we prove an analog of Theorem A in equal characteristic $p>0$
in which the infinite set of points consisting of all $T$ as above is replaced by the finite set
$\{x_i\}_{i=1}^m$.
%Further, we prove the following (cf. Proposition 3.4 in equal characteristic $p>0$ and Proposition 3.6 in mixed characteristics).
%\proclaim {Theorem B} Assume $p=\char(k)\ge 0$, $\ell$ is a prime integer (possibly equal to $p$ if $p>0$), and with no restrictions on $\char (K)$. 
%Then the pro-$\ell$ group $\pi_1(\Cal U\setminus S)^{\geo,\ell}$ is free 
%of infinite (resp. finite) rank if $\ell=p> 0$ (if $\ell\neq p$). 
%\endproclaim
Further, we prove the following (cf. Theorem 3.7) which, in case $\char(k)=0$, gives a description of the structure of (the full) $\pi_1(\Cal U\setminus S)^{\geo}$ in the equal characteristic $0$ case.

\proclaim {Theorem B} Assume $\char (k)=p\ge 0$ with no restriction on $\char (K)$. Then the 
morphism $\Cal U\to X_K$ induces (via the rigid GAGA functor) a continuous homomorphism 
$\pi_1(\Cal U\setminus S)^{\geo,p'}\to \pi_1(X_K\setminus (\{x_i\}_{i=1}^m\cup S))^{\geo,p'}$ between the maximal prime-to-$p$ quotients of $\pi_1(\Cal U\setminus S)^{\geo}$
and $\pi_1(X_K\setminus (\{x_i\}_{i=1}^m\cup S))^{\geo}$; respectively, which is an isomorphism.
In particular, if $S(\overline K)=\{y_1,\cdots,y_r\}$ has cardinality $r\ge 0$ then $\pi_1(\Cal U\setminus S)^{\geo,p'}$ is (pro-)prime-to-$p$ free 
on $2g+m+r-1$ generators
and can be generated by $2g+m+r$ generators $\{a_1,\cdots,a_g,b_1,\cdots,b_g,\sigma_1,\cdots,\sigma_m,
\tau_1,\cdots,\tau_r\}$
subject to the unique relation $\prod_{j=1}^g[a_j,b_j]\prod _{i=1}^m\sigma _i\prod_{t=1}^r \tau _t=1$, where $\sigma_i$ (resp $\tau _t$) is a generator of inertia at $x _i$
(resp. $y_t$) and $g\defeq g_{X_K}$ is the arithmetic genus of $X_K$ (also called the genus of the affinoid $\Cal U$). 
\endproclaim

Note that unless $\char(k)=0$ the profinite group $\pi_1(\Cal U\setminus S)^{\geo}$ is {\it not free} (neither is it finitely generated) as the ranks of its maximal pro-$\ell$ quotients 
can be different for different primes $\ell$ (cf. Theorem A, and its analog Theorem 3.3 in equal characteristic $p>0$).  In this sense Theorem A (and Theorem 3.3) is an optimal result one can prove regarding the structure of the {\it full} geometric 
fundamental group of a $p$-adic smooth affinoid curve. Also there is no analog to Theorem A if $\char(k)=p>0$, for the full $\pi_1^{\geo}$, 
where one replaces the infinite union of the finite sets of points $T$ 
(as in the statement of Theorem A) by a single fixed finite set of points $\widetilde T\subset \bigcup_{i=1}^m \Cal D_{i}^o$ (cf. Remark 3.8(ii)).

Next, we outline the content of the paper. In $\S1$ we collect some well-known background material.
In $\S2$ we explain how one defines the \'etale fundamental group of a rigid analytic $K$-affinoid space (cf. 2.1), and recall the rigid analog of Runge's theorem proven by Raynaud (cf. 2.2).
We then investigate in 2.3 the structure of a certain quotient of the geometric \'etale fundamental group of an annulus of thickness $0$. In $\S3$ we 
investigate the structure of the geometric \'etale fundamental group of a smooth affinoid $p$-adic curve and prove Theorem A (as well as its analog Theorem 3.3 if $\char(K)=p>0$), and Theorem B.

In [Garuti] Garuti investigated, among others, the structure of the pro-$p$ geometric fundamental group of a rigid closed $p$-adic annulus of thickness $0$ and proved an analogue of
Theorem A in this special case.

\subhead
Acknowledgment
\endsubhead
I would like to thank the referee for his/her careful reading of the paper and the many suggestions which helped improve the presentation of the paper.

\subhead
Notations
\endsubhead
In this paper $K$ is a complete discrete valuation ring, $R$ its valuation ring, $\pi$ a uniformising parameter, 
and $k\defeq R/\pi R$ the residue field of characteristic $p\ge 0$ which we assume to be algebraically closed. 

We refer to [Raynaud], $3$, for the terminology we will use 
concerning $K$-rigid analytic spaces, $R$-formal schemes, as well as the link between formal and rigid geometry. 
For an $R$-(formal) scheme $X$ we will denote by
$X_K\defeq X\times _RK$ (resp. $X_k\defeq X\times _Rk$) the generic (resp. special) fibre of $X$ 
(the generic fibre is understood in the rigid analytic sense in the case where $X$ is a formal scheme).
Moreover, if $X=\Spf A$ 
is an affine formal $R$-scheme of finite type we denote by $X_K\defeq \Sp (A\otimes _RK)$ the associated $K$-rigid affinoid space and
will also denote, when there is no risk of confusion, by $X_K$
the affine scheme $X_K\defeq \Spec (A\otimes _RK)$.

A formal (resp. algebraic) $R$-curve is an $R$-formal scheme of finite type (resp. scheme of finite type)
flat and separated whose special fibre is equidimensional of dimension $1$.
For a $K$-scheme (resp. $K$-rigid analytic space) $X$ and $L/K$ a field extension (resp. a finite extension) we write
$X_L\defeq X\times _KL$ which is an $L$-scheme (resp. an $L$-rigid analytic space). If $X$ is a proper and normal
formal $R$-curve we also denote, when there is no risk of confusion, by $X$ the algebraisation of $X$ which is an algebraic normal and proper $R$-curve and by $X_K$ 
the proper and normal algebraic $K$-curve associated to the rigid $K$-curve $X_K$ via the rigid GAGA functor.

For a profinite group $H$ and a prime integer $\ell$ we denote by $H^{\ell}$ the maximal pro-$\ell$ quotient of $H$,
and $H^{\ell'}$ the maximal prime-to-$\ell$ quotient of $H$.

All scheme cohomology groups $H_{\et}^1(\ ,\Bbb Z/\ell\Bbb Z)$ in this paper are \'etale cohomology groups.

\subhead
\S 1 Background
\endsubhead

\subhead
1.1
\endsubhead
Let $p>1$ be a prime integer. We recall some well-known facts 
on profinite pro-$p$ groups.
First, we recall the following characterisations of free pro-$p$ groups.

\proclaim {Proposition 1.1} Let $G$ be a profinite pro-$p$ group. Then the following properties are equivalent.

\noindent
(i)\ $G$ is a free pro-$p$ group.

\noindent
(ii)\ The $p$-cohomological dimension of $G$ satisfies $\cd_p(G)\le 1$.

\noindent
In particular, a closed subgroup of a free pro-$p$ group is free.

\endproclaim

\demo{Proof} Well-known (cf. [Serre], and [Ribes-Zalesskii], Theorem 7.7.4).
\qed
\enddemo

Next, we recall the notion of a {\it direct factor} of a free pro-$p$ group 
(cf. [Garuti], 1, the discussion preceding  Proposition 1.8).

\proclaim {Definition/Lemma 1.2 (Direct factors of free  pro-$p$ groups)} 
Let $F$ be a free pro-$p$ group, $H\subseteq F$ a closed 
subgroup, and $\iota: H\to F$ the natural homomorphism. We say that 
$H$ is a direct factor of $F$ if there exists a continuous homomorphism 
$s:F\to H$ such that $s\circ \iota=\id_H$ ($s$ is necessarily surjective). 
There exists then a (non unique) closed subgroup $N$ of $F$ such that $F$ is isomorphic to the free direct product $H\star N$.
We will refer to such a subgroup $N$ as a supplement of $H$.
\endproclaim

\demo{Proof} In what follows we consider $\Bbb Z/p\Bbb Z$ as the trivial discrete module. 
Let $s:F\to H$ be a left inverse to $\iota$ as above which induces a retraction $h^1(s):H^1(H,\Bbb Z/p\Bbb Z)\to H^1(F,\Bbb Z/p\Bbb Z)$ of the map
$h^1(\iota):H^1(F,\Bbb Z/p\Bbb Z)\to H^1(H,\Bbb Z/p\Bbb Z)$ induced by $\iota$; in particular $h^1(s)$ is injective. Let $M$ be a supplement of (the image via $h^1(s)$ of) 
$H^1(H,\Bbb Z/p\Bbb Z)$ in $H^1(F,\Bbb Z/p\Bbb Z)$ and $M^{\wedge}$ the corresponding subgroup of $F/F^{\star}$ where $F^{\star}$ is the Frattini 
subgroup of $F$ (recall that $F/F^{\star}$ is the Pontrjagin 
dual of $H^1(F,\Bbb Z/p\Bbb Z)$). Let $\{\tilde g_i\}_i$ be a minimal set of generators of $M^{\wedge}$ and $\{g_i\}_i\subset F$ a lift of the $\{\tilde g_i\}_i$.
Write $N$ for the closed subgroup of $F$ generated by the $\{g_i\}_i$ (which is free pro-$p$) and $H\star N$ the free product of $H$ and $N$. Then the natural morphism 
$H\star N\to F$ is an isomorphism as it is an isomorphism on the cohomology with coefficients in $\Bbb Z/p\Bbb Z$ (cf. [Ribes-Zalesskii], Proposition 7.7.2).
\qed
\enddemo

Note that in the discussion before Proposition 1.8 in [Garuti], and with the notations in Definition/Lemma 1.2,
it is stated that $F$ is isomorphic to the free product $H\star \Ker(s)$. This is not necessarily the case as the induced map 
on cohomology is not necessarily an isomorphism as claimed in loc. cit. A similar inaccurate statement occurs in the proof of Proposition 1.8 of loc. cit., but this doesn't affect the 
validity of this Proposition.

One has the following cohomological characterisation of direct factors of free 
pro-$p$ groups.

\proclaim {Proposition 1.3}
Let $H$ be a pro-$p$ group, $F$ a free pro-$p$ group, and 
$\sigma :H\to F$ a continuous homomorphism. Assume that the map induced by $\sigma$ on 
cohomology
$$h^1(\sigma):H^1(F,\Bbb Z/p\Bbb Z)\to H^1(H,\Bbb Z/p\Bbb Z)$$
is surjective, where
$\Bbb Z/p\Bbb Z$ is considered as a trivial discrete module. 
Then $\sigma$ induces an isomorphism $H\isom \sigma(H)$ and $\sigma(H)$ is a direct factor of $F$. 
We say that $\sigma$ makes $H$ into a direct factor of $F$.
\endproclaim

\demo{Proof} cf. [Garuti], Proposition 1.8.
\qed
\enddemo

Next, we consider the notion of a {\it semi-direct factor} of a profinite group.

\definition {Definition 1.4 (Semi-direct factors of profinite groups)} Let $G$ be a profinite group, $H\subseteq G$ a closed subgroup,
and $\iota: H\to G$ the natural homomorphism. We say that 
$H$ is a {\it semi-direct factor} of $G$ if there exists a continuous homomorphism 
$s:G\to H$ such that $s\circ \iota=\id_H$ ($s$ is necessarily surjective). 
\enddefinition

Note that $G$ is the semi-direct product of $\Ker s$ and $H$. Unlike the pro-$p$ case, a semi-direct factor of a free profinite group is not necessarily free.

\proclaim{Lemma 1.5} Let $\tau:H\to G$ be a continuous homomorphism between profinite groups. Write 
$H=\underset{j\in J} \to {\varprojlim}\ H_j$ as the projective limit of the inverse system $\{H_j,\phi_{j'j},J\}$
of finite quotients $H_j$ of $H$ with index set $J$.
Suppose there exists,  $\forall j\in J$, a surjective homomorphism $\psi_j :G\twoheadrightarrow H_j$ such that 
$\psi_j \circ \tau:H\twoheadrightarrow H_j$ is the natural map and
$\psi_j=\phi _{j'j}\circ \psi _{j'}$ whenever this makes sense.
Then $\tau$ induces an isomorphism $H\isom \tau(H)$ and $\tau(H)$ is a semi-direct factor of $G$. 
We say that $\tau$ makes $H$ into a semi-direct factor of $G$.
\endproclaim

\demo{Proof} Indeed, the $\{\psi_j\}_{j\in J}$ give rise to a continuous (necessarily surjective) homomorphism $\psi :G\to H$ which is a left 
inverse of $\tau$.
\qed
\enddemo

\subhead
\S 2. Geometric fundamental groups of annuli of thickness zero
\endsubhead
In this section we explain how one defines the \'etale fundamental group of a rigid analytic $K$-affinoid space 
(cf. 2.1), and recall the rigid analog of Runge's Theorem proven by Raynaud (cf. 2.2).
We then investigate, in 2.3, the structure of a certain quotient of the geometric \'etale fundmantal group of an annulus of thickeness $0$.
The main results in this section are inspired from [Garuti], $\S2$.

\subhead
2.1
\endsubhead
First, we explain how one defines the \'etale fundamental group of a rigid analytic $K$-affinoid space.
Let  $U=\Spf A$ 
be an affine $R$-formal scheme which is topologically of finite type. Thus, 
$A$ is a $\pi$-adically complete noetherian $R$-algebra. Let 
$\Cal A\defeq A\otimes _RK$ be the corresponding Tate algebra and 
$\Cal U\defeq \Sp \Cal A$ 
the associated  $K$-rigid analytic affinoid space which is the generic fibre of 
$U$ in the sense of Raynaud (cf. [Raynaud], 3). 
Assume that the affine scheme $\Spec \Cal A$ is (geometrically) normal and 
geometrically connected. Let $\eta$ be a geometric point of  
$\Spec \Cal A$ above its generic point. Then $\eta$ 
determines an algebraic closure $\overline K$ of $K$ and a geometric 
point of $\Spec (\Cal A\times _K\overline K)$ which we will also denote $\eta$.

\definition {Definition 2.1.1 (\'Etale Fundamental Groups of Affinoid Spaces)} 
(See also [Garuti], D\'efinition 2.2 and D\'efinition 2.3).
We define the {\it \'etale fundamental group} of $\Cal U$ with base point $\eta$ by
$$\pi_1(\Cal U,\eta)\defeq \pi_1(\Spec \Cal A,\eta),$$
where $\pi_1(\Spec \Cal A,\eta)$ is the \'etale fundamental group of the connected 
scheme $\Spec \Cal A$ with base point $\eta$ in the sense of Grothendieck (cf. [Grothendieck], V).
Thus, $\pi_1(\Cal U,\eta)$ classifies finite coverings 
$\Spec {\Cal B}\to \Spec {\Cal A}$ where $\Cal B$ is a finite \'etale $\Cal A$-algebra.
There exists a continuous surjective homomorphism
$\pi_1(\Cal U,\eta)\twoheadrightarrow \Gal (\overline K/K).$
We define the {\it geometric \'etale fundamental group} $\pi_1(\Cal U,\eta)^{\geo}$ of $\Cal U$ so 
that the following sequence is exact
$$1\to \pi_1(\Cal U,\eta)^{\geo}\to \pi_1(\Cal U,\eta)\to \Gal (\overline K/K)\to 1.$$
\enddefinition

\definition {Remark 2.1.2} If $L/K$ is a finite field extension contained in 
$\overline K/K$, and $\Cal U_L\defeq \Cal U\times _KL$ 
is the affinoid $L$-rigid analytic 
space obtained from $\Cal U$ by extending scalars, then we have a 
commutative diagram:
$$
\CD 
1 @>>>  \pi _1(\Cal U_L,\eta)^{\geo} @>>> \pi_1(\Cal U_L,\eta) @>>> \Gal (\overline K/L) 
@>>> 1 \\
@.  @VVV   @VVV     @VVV \\
1 @>>> \pi _1 (\Cal U,\eta)^{\geo} @>>> \pi_1(\Cal U,\eta) @>>> \Gal (\overline K/K) 
@>>> 1\\
\endCD
$$ 
where the two right vertical maps are injective homomorphisms and the left 
vertical map is an isomorphism. 
The geometric fundamental group $\pi_1(\Cal U,\eta)
^{\geo}$ is strictly speaking {\it not} the fundamental group of a rigid analytic 
space (since $\overline K$ is not complete). It is, however, the projective limit 
of fundamental groups of rigid affinoid spaces. More precisely, there exists 
an isomorphism
$$\pi_1(\Cal U,\eta)^{\geo} \isom \underset{L/K} \to {\varprojlim}\pi_1(\Cal U\times _K L,\eta),$$
where the limit is taken over all finite extensions $L/K$ contained in 
$\overline K$.
\enddefinition

Similarly, if $\Cal U$ above is a geometrically connected and (geometrically) normal {\it affinoid $K$-curve}, 
and $S$ is a finite set of points of $\Cal U$ (cf. [Raynaud], 3.1, for the definition of points of a rigid analytic space), 
we define the \'etale fundamental group $\pi_1(\Cal U\setminus S,\eta)$
of $\Cal U\setminus S$ with base point $\eta$ which is a profinite group and classifies finite coverings $\Spec {\Cal B}\to \Spec {\Cal A}$, where $\Cal B$ 
is a finite $\Cal A$-algebra which is \'etale above $\Cal U\setminus S$. In this case we have an exact sequence
$$1\to \pi_1(\Cal U\setminus S,\eta)^{\geo}\to \pi_1(\Cal U\setminus S,\eta)\to \Gal (\overline K/K)\to 1,$$
where $\pi_1(\Cal U\setminus S,\eta)^{\geo}\defeq  \Ker \lgroup \pi_1(\Cal U\setminus S,\eta)\twoheadrightarrow \Gal (\overline K/K)\rgroup$, and a similar description of $\pi_1(\Cal U\setminus S,\eta)^{\geo}$ to that of $\pi_1(\Cal U,\eta)^{\geo}$ given in Remark 2.1.2.

\subhead
2.2
\endsubhead
Next, we recall the rigid analog of Runge's Theorem proven by Raynaud. Let $X_K$ be a proper, smooth, and geometrically connected algebraic $K$-curve. We denote by $X_K^{\rig}$ the associated $K$-rigid analytic proper and smooth curve.
Let $\Cal U\hookrightarrow X_K^{\rig}$ be an open affinoid subspace of $X_K^{\rig}$ (cf. [Raynaud], 3.1). The following is well-known (cf. [Raynaud], Proposition 3.5.1).

\proclaim {Proposition 2.2.1}
The complement $\Cal W\defeq X_K^{\rig}\setminus \Cal U$ has a natural structure of
an (non quasi-compact) open rigid subspace of $X_K^{\rig}$ which is an increasing union of open quasi-compact rigid subspaces of $X_K^{\rig}$.
The rigid space $\Cal W$ has a finite number of connected components $\{\Cal W_i\}_{i\in I}$. For each $i\in I$, let $x_i\in \Cal W_i$ be a point 
(in the sense of [Raynaud], 3.1) and write $U_K\defeq X_K\setminus \{x_i\}_{i\in I}$
which is an affine $K$-curve. Then there exists a canonical affine and normal $R$-scheme $U^{\alg}$  of finite type 
such that $(U^{\alg})_K=U_K$, and if $\widetilde U=\Spf A$ denotes the formal completion of $U^{\alg}$ for the $\pi$-adic topology then the  
generic fibre $\widetilde U_K=\Sp \Cal A$ of $\widetilde U$ (in the sense of [Raynaud], 3.1), where $\Cal A\defeq A\otimes _RK$, is the rigid affinoid $K$-curve $\Cal U$.
\endproclaim

As a consequence one obtains the following version of Runge's Theorem for rigid $K$-curves (cf. [Raynaud], Corollaire 3.5.2).

\proclaim{Proposition 2.2.2 (Runge's Theorem)} We use the same notations as in Proposition 2.2.1. Then 
the ring of regular functions on the affine curve $U_K$ has a dense image in the ring of holomorphic functions on $\Cal U$.
More generally, a coherent sheaf $M_K$ on $U_K$ induces a coherent sheaf $\Cal M$ on $\Cal U$ and the image of the sections of $M_K$ on $U_K$ is dense in the space 
of sections of $\Cal M$ on $\Cal U$.
\endproclaim

We will refer to a pair $(\Cal U,U_K)$ as in Proposition 2.2.1 as a {\it Runge pair}.

\subhead
2.3
\endsubhead
In this section we investigate the structure of a certain quotient of the geometric \'etale fundamental group of an annulus of thickness $0$.
Let $D=\Spf R<Z>$ be the formal standard closed disc and $\Cal D\defeq D_K=\Sp K<Z>$ its generic fibre which is the standard closed rigid analytic disc centred at the point $"Z=0"$. 
Given an integer $n\ge 0$ consider the formal closed disc $D_n\defeq \Spf \frac {R<Z,Y>}{(Z-\pi^nY)}$ and its generic fibre $\Cal D_n\defeq D_{n,K}=\Sp \frac {K<Z,Y>}{(Z-\pi^nY)}$
(recall $\pi$ is a uniformiser of $R$). 
The natural embedding $\Cal D_n\subset \Cal D_0=\Cal D$ induces an identification between the points of $\Cal D_n$ and the closed disc $\{x\in \Cal D,\ \vert Z(x) \vert \le \vert \pi \vert^n\}$.
We also consider the formal annulus $C_n\defeq \Spf \frac {R<Z,Y,W>}{(Z-\pi^nY,YW-1)}$ and its generic fibre $\Cal C_n\defeq C_{n,K}=\Sp \frac {K<Z,Y,W>}{(Z-\pi^nY,YW-1)}$. 
The natural embedding $\Cal C_n\subset \Cal D_n$ induces an identification between the points of $\Cal C_n$ and the closed annulus of thickeness zero $\{x\in \Cal D,\ \vert Z(x) 
\vert = \vert \pi \vert^n\}$. 

Let $U_K\defeq \Bbb G_{m,K}=\Spec \frac {K[Z,V]}{(ZV-1)}$ and $X_K=\Bbb P^1_K$ its smooth compactification, with function field $K(Z)$. 
We have natural embeddings $\Cal C_n\subset \Cal D_n\subset \Bbb (X_K)^{\rig}$. (Here we consider the rigid analytic structure on $X_K$ 
arising from the admissible covering  $\{x\in \Bbb P^1_K,\ \vert Z(x) 
\vert \le \vert \pi \vert^n\}\bigcup \{x\in \Bbb P^1_K,\ \vert Z(x)\vert \ge \vert \pi \vert^n\}$.)
Let $\eta$ be a geometric point of $\Cal C_n$ as in 2.1 which induces a geometric point of $\Cal D_n$, $U_K$, $X_K$ and $X_{\overline K} \defeq X_K\times _K {\overline K}$ 
(which we also denote $\eta$). There exist continuous homomorphisms $\phi_n:\pi_1(\Cal C_n,\eta)\to \pi_1(\Cal D_n\setminus \{0\},\eta)$ and $\psi_n:\pi_1
(\Cal D_n\setminus \{0\},\eta)
\to \pi_1(U_{K}, \eta)$ (via the rigid GAGA functor) which induce continuous homomorphisms 
$\phi_n^{\geo}:\pi_1(\Cal C_n,\eta)^{\geo}\to \pi_1(\Cal D_n\setminus \{0\},\eta)^{\geo}$ and $\psi_n^{\geo}:
\pi_1(\Cal D_n\setminus \{0\},\eta)^{\geo}\to \pi_1(U_{K}, \eta)^{\geo}$.

\proclaim {Proposition 2.3.1} Let $p=\char(k)\ge0$ with no restriction on $\char(K)$. Then the homomorphisms
$\phi_n^{\geo,p'}:\pi_1(\Cal C_n,\eta)^{\geo,p'}\to  \pi_1(\Cal D_n\setminus \{0\},\eta)^{\geo,p'}$ and $\psi_n^{\geo,p'}:\pi_1(\Cal D_n\setminus \{0\},\eta)^{\geo,p'}
\to \pi_1(U_{K}, \eta)^{\geo,p'}$ (induced by  $\phi_n^{\geo}$ and $\psi_n^{\geo}$; respectively)
are isomorphisms. In particular, both $\Gamma\defeq \pi_1(\Cal C_n,\eta)^{\geo,p'}$ and $\widetilde \Gamma \defeq \pi_1(\Cal D_n\setminus \{0\},\eta)^{\geo,p'}$ 
are isomorphic to the maximal prime-to-$p$ quotient $\hat \Bbb Z^{p'}$ of $\hat \Bbb Z$.
\endproclaim

\demo{Proof} The last assertion follows from the first, and the well-known fact (since $U_K=\Bbb G_{m,K}$) that $\pi_1(U_{K}, \eta)^{\geo,p'}$ is 
isomorphic to $\hat \Bbb Z^{p'}$. 
%To prove the first assertion we can assume, without loss of generality, that $n=0$.
Let $C_{n,k}\defeq \Spec \frac {k[y,w]}{(yw-1)}=\Bbb G_{m,k}$ be the special fibre of 
$C_n$ and $\beta$ a geometric point of $C_{n,k}$ which induces a geometric point of $C_n$ noted also $\beta$. There exist continuous homomorphisms
$\pi_1(\Cal C_n,\eta)^{\geo}\to \pi_1(C_n,\beta)^{\geo}\to \pi_1(C_{n,k},\beta)^{\geo}$, where the first map is surjective (a geometrically connected 
\'etale cover of $C_n$ induces by passing to the generic fibres a geometrically connected \'etale cover of $\Cal C_n$), and the second map is 
the inverse of the natural map $\pi_1(C_{n,k},\beta)^{\geo}\to \pi_1(C_n,\beta)^{\geo}$ which is an isomorphism (cf. [SGA1], Expos\'e I, Corollaire 8.4).
The composite map $\pi_1(\Cal C_n,\eta)^{\geo}\twoheadrightarrow \pi_1(C_{n,k},\beta)^{\geo}$ is a surjective specialisation homomorphism,
which induces a surjective specialisation homomorphism $\pi_1(\Cal C_n,\eta)^{\geo,p'}\twoheadrightarrow \pi_1(C_{n,k},\beta)^{\geo,p'}$ between the respective
prime-to-$p$ parts. One can show, using Abhyankar's lemma (cf. loc. cit. Expos\'e X, Lemma 3.6)
and the theorem of purity of Zariski (cf. loc. cit. Expos\'e X, Th\'eor\`eme 3.1), that this latter map is an isomorphism (similar arguments used in loc. cit., 
Expos\'e X, in order to prove Th\'eor\`eme 3.8 and Corollaire 3.9). 
On the other hand $\pi_1(C_{n,k},\beta)^{\geo,p'}\isom \pi_1(\Bbb G_{m,k})^{\geo,p'}\isom \hat \Bbb Z^{p'}$. 
Thus, $\pi_1(\Cal C_n,\eta)^{\geo,p'}\isom \hat \Bbb Z^{p'}$.

Further, the composite homomorphism 
$\pi_1(\Cal C_n,\eta)^{\geo,p'}\to \pi_1(\Cal D_n\setminus \{0\},\eta)^{\geo,p'}\to \pi_1(U_{K}, \eta)^{\geo,p'}$ is an isomorphism. 
Indeed, this composite map 
is surjective since
a finite Galois \'etale cover $V_K\to U_K$ of order prime-to-$p$, with $V_K$ geometrically connected, extends to a finite (cyclic) Galois cover 
$Z\to X_K$ (totally) ramified above $0$ and $\infty$ and the corresponding Galois cover $Z^{\rig}\to (X_K)^{\rig}$ of rigid curves restricts 
(after possibly a finite extension of $K$) to an \'etale cover $\Cal V_n\to \Cal C_n$ with $\Cal V_n$ geometrically connected. More precisely,
the rigid cover  $Z^{\rig}\to (X_K)^{\rig}$ induces in reduction, via suitable formal models, a finite \'etale cover  $V_k\to C_{n,k}$ with $V_k$ geometrically connected,
hence $\Cal V_n$ is geometrically connected. 
As both $\pi_1(\Cal C_n,\eta)^{\geo,p'}$ and $\pi_1(U_{K}, \eta)^{\geo,p'}$ are isomorphic to $\hat \Bbb Z^{p'}$, the composite surjective map 
$\pi_1(\Cal C_n,\eta)^{\geo,p'}\to \pi_1(U_{K}, \eta)^{\geo,p'}$ is an isomorphism. 
Further, the map  
$\pi_1(\Cal C_n,\eta)^{\geo,p'}\to \pi_1(\Cal D_n\setminus \{0\},\eta)^{\geo,p'}$
is surjective (same argument as the one used above for the surjectivity of the map $\pi_1(\Cal C_n,\eta)^{\geo,p'}\to \pi_1(U_{K}, \eta)^{\geo,p'}$).
We then deduce that the maps 
$\pi_1(\Cal C_n,\eta)^{\geo,p'}\to \pi_1(\Cal D_n\setminus \{0\},\eta)^{\geo,p'}$ and
$\pi_1(\Cal D_n\setminus \{0\},\eta)^{\geo,p'}\to \pi_1(U_{K}, \eta)^{\geo,p'}$ are isomorphisms as claimed.
\qed
\enddemo

\proclaim {Proposition 2.3.2} Assume $\char (K)=p>0$. Then the homomorphism
$\psi_n^{\geo,p}:\pi_1(\Cal D_n\setminus \{0\},\eta)^{\geo,p}
\to \pi_1(U_{K}, \eta)^{\geo,p}$ (induced by  $\psi_n^{\geo}$) makes $\pi_1(\Cal D_n\setminus \{0\},\eta)^{\geo,p}$ 
into a direct factor of $\pi_1(U_{K}, \eta)^{\geo,p}$ and $\pi_1(\Cal D_n\setminus \{0\},\eta)^{\geo,p}$ is a free pro-$p$ group of infinite rank.
Furthermore, the homomorphism
$\phi_n^{\geo,p}:\pi_1(\Cal C_n,\eta)^{\geo,p}\to  \pi_1(\Cal D_n\setminus \{0\},\eta)^{\geo,p}$ (induced by  $\phi_n^{\geo}$) makes
$\pi_1(\Cal C_n,\eta)^{\geo,p}$ into a direct factor of $\pi_1(\Cal D_n\setminus \{0\},\eta)^{\geo,p}$ and
$\pi_1(\Cal C_n,\eta)^{\geo,p}$ is a free pro-$p$ group of infinite rank.
\endproclaim

\demo{Proof} First, note that $\pi_1(U_{K}, \eta)^{\geo,p}$ is free
since $U_K$ is an affine scheme of characteristic $p>0$ (cf. [Serre1], Proposition 1).
Next, we prove the first assertion.
Using Proposition 1.3, we need to show that the map 
$H^1(\pi_1(U_{K}, \eta)^{\geo},\Bbb Z/p\Bbb Z)\to H^1(\pi_1(\Cal D_n\setminus 0,\eta)^{\geo},\Bbb Z/p\Bbb Z)$
induced by $\psi_n^{\geo}$ on cohomology is surjective. 
Let $f:\Cal Z\to \Cal D_n$ be a generically $\Bbb Z/p\Bbb Z$-torsor which is \'etale outside $0$ with $\Cal Z$ geometrically connected.
Let $n'>n$ an integer and $\Cal X\defeq \Cal D_n\setminus \Cal D_{n'}^o$; where $\Cal D_{n'}^o\defeq \Cal D_{n'}\setminus \Cal C_{n'}$ is an open disc, 
which is an affinoid subdomain of $\Cal D_n$.
Let $\tilde f:\Cal Y\to \Cal X$ be the restriction of $f$ which is an \'etale $\Bbb Z/p\Bbb Z$-torsor.
For $n'>>n$, $f^{-1}(\Cal D_{n'}^o)$ is geometrically connected since $f$ is totally ramified 
above $0$, and $\Cal Y$ is then geometrically connected, which we will assume from now on.
(More precisely, the pre-image of $0\in \Cal D_n(K)$ in $\Cal Z$ consists of a single point $z\in \Cal Z(K)$ as $f$ is totally ramified above $0$.
By passing to a formal model of $\Cal Z$, its minimal desingularisation, and the quotient of the latter by the action of the Galois group $\Bbb Z/p\Bbb Z$ 
of the covering $f$, one sees that $f^{-1}(\Cal D_{n'}^o)$ is an open disc for $n'>>n$.)
By Artin-Schreier theory the torsor $\tilde f$ is given by an Artin-Schreier equation $\alpha ^p-\alpha=g$ where $g$ is a holomorphic function on $\Cal X$.
The pair $(\Cal X,U_K)$ is a Runge pair.
The function $g$ can be approximated by a regular function $\tilde g$ on $U_K$ (cf. Proposition 2.2.2). For $\tilde g$ close to $g$ the equation 
$\alpha^p-\alpha=\tilde g$ 
defines a $\Bbb Z/p\Bbb Z$-\'etale torsor  $f':Z_K\to U_K$ whose pull-back to $\Cal X$ is isomorphic to $\tilde f$.
In particular, $Z_K$ is geometrically connected.
More precisely, for $\tilde g$ close to $g$ (for example if $\vert \vert \tilde g-g\vert \vert <1$, where $\vert\vert \ \ \vert\vert$ is the supremum norm on 
$\Cal X$) then $\sum_{t\ge 0} (\tilde g-g)^{p^t}$ converges to a holomorphic function $h$ on $\Cal X$ and $g=h^p-h+\tilde g$ in $\Cal X$. Hence the class of $f'$ in
$H^1(\pi_1(U_{K}, \eta)^{\geo},\Bbb Z/p\Bbb Z)$ maps to the class of $f$ in $H^1(\pi_1(\Cal D_n\setminus 0,\eta)^{\geo},\Bbb Z/p\Bbb Z)$.
The assertion that $\pi_1(\Cal D_n\setminus \{0\},\eta)^{\geo,p}$ is free follows from loc. cit. (cf. Proposition 1.1).
Similarly, using the fact that $(\Cal C_n,U_K)$ is a Runge pair, one proves that the natural map 
$H^1(\pi_1(U_{K}, \eta)^{\geo},\Bbb Z/p\Bbb Z)\to H^1(\pi_1(\Cal C_n,\eta)^{\geo},\Bbb Z/p\Bbb Z)$
is surjective, hence the second assertion follows from Proposition 1.3. Finally, the assertions on infinite rank follow from the facts that
$H^1(\pi_1(\Cal D_n\setminus 0,\eta)^{\geo},\Bbb Z/p\Bbb Z)$ and $H^1(\pi_1(\Cal C_n,\eta)^{\geo},\Bbb Z/p\Bbb Z)$ are infinite dimensional $\Bbb F_p$-vector spaces.
\qed
\enddemo

Next, we investigate the structure of the maximal pro-$p$ quotients of $\pi_1(\Cal C_n,\eta)^{\geo}$ and $\pi_1(\Cal D_n\setminus T,\eta)^{\geo}$ ( 
$T\subset \Cal D_n\setminus \Cal C_n$ is a finite set of points) in the mixed characteristic case.
Let $T$ be a finite set of points of $\Cal D_n\setminus \Cal C_n$ and $S$ a finite set of points of $X_K^{\rig}\setminus \Cal D_n$.
We view $T\cup S \subset X_K$ as a closed subscheme of $X_K$ and write 
$(T\cup S)_L\defeq (T\cup S)\times _KL$ if $L/K$ is a sub-extension of $\overline K/K$.
We also denote by $\pi_1(X_L\setminus (T\cup S)_L, \eta)$ the \'etale fundamental group 
of  $X_L\setminus (T\cup S)_L$ with base point $\eta$. The natural embedding 
$\Cal D_{n,L}\defeq \Cal D_n\times _KL\to X_L^{\rig}$ 
induces (via the rigid GAGA functor) a continuous homomorphism
$\pi_1(\Cal D_{n,L}\setminus T_L,\eta)\to \pi_1(X_L\setminus (T\cup S)_L, \eta)$,
and by passing to the projective limit a homomorphism (cf. Remark 2.1.2)
$$\pi_1(\Cal D_n\setminus T,\eta)^{\geo}\to \pi_1(X_{\overline K}\setminus (T\cup S)_{\overline K}, \eta) 
\defeq \underset{L/K} \to {\varprojlim} \pi_1(X_L\setminus (T\cup S)_L, \eta),$$
where $L/K$ runs over all finite extensions contained in $\overline K$.

Let $\ell$ be a prime integer.
The above homomorphism 
$\pi_1(\Cal D_n\setminus T,\eta)^{\geo}\to \pi_1(X_{\overline K}\setminus (T\cup S)_{\overline K}, \eta)$
 induces homomorphisms
$ \pi_1(\Cal D_n\setminus T,\eta)^{\geo,\ell}\to \pi_1(X_{\overline K}\setminus (T\cup S)_{\overline K}, \eta)^{\ell}$ 
and (by passing to the projective limit over all $S$ as above)
$$\phi_{n,T}:\pi_1(\Cal D_n\setminus T,\eta)^{\geo,\ell}\to \underset {S}\to {\varprojlim}\ \pi_1(X_{\overline K}\setminus (T\cup S)_{\overline K}, \eta)^{\ell},$$ 
which induces a homomorphism 
$$\phi_n\defeq \underset{T} \to {\varprojlim}\ \phi_{n,T}: \underset{T}\to {\varprojlim}\ \pi_1(\Cal D_n\setminus T,\eta)^{\geo,\ell}\to 
\underset{(T,S)}\to {\varprojlim}\ \pi_1(X_{\overline K}\setminus (T\cup S)_{\overline K},\eta)^{\ell},$$ 
where the limit is taken over all finite sets of points $T$ and $S$ as above
and
$\underset{(T,S)}\to {\varprojlim}\ \pi_1(X_{\overline K}\setminus (T\cup S)_{\overline K},\eta)^{\ell}\defeq 
\underset{T}\to {\varprojlim}\ \left(\underset{S}\to {\varprojlim}\ \pi_1(X_{\overline K}\setminus (T\cup S)_{\overline K},\eta)^{\ell}\right)$. 
The profinite groups $\underset {S}\to {\varprojlim}\ \pi_1(X_{\overline K}\setminus (T\cup S)_{\overline K}, \eta)^{\ell}$ and 
$\underset{(T,S)}\to {\varprojlim}\ \pi_1(X_{\overline K}\setminus (T\cup S)_{\overline K},\eta)^{\ell}$
are {\it free} pro-$\ell$ groups if $\ell\neq \char (K)$ (as follows from [Grothendieck], Expos\'e XIII, Corollaire 2.12). 
Note that for each finite set $T\subset \Cal D_n\setminus \Cal C_n$ as above we have a continuous homomorphism
$\pi_1(\Cal C_n,\eta)^{\geo}\to \pi_1(\Cal D_n\setminus T,\eta)^{\geo}$ which induces homomorphisms 
$\psi_{n,T}:\pi_1(\Cal C_n,\eta)^{\geo,\ell}\to \pi_1(\Cal D_n\setminus T,\eta)^{\geo,\ell}$ and (by passing to the projective limit over all $T$ as above)
$$\psi_{n}\defeq \underset{T} \to {\varprojlim}\ \psi_{n,T}:\pi_1(\Cal C_n,\eta)^{\geo,\ell}\to \underset {T}\to {\varprojlim}\ \pi_1(\Cal D_n\setminus T,\eta)^{\geo,\ell}.$$

\proclaim {Proposition 2.3.3} Assume $\char(K)=0$ and $\char(k)=p>0$.
Then the continuous homomorphism $\phi_{n,T}:\pi_1(\Cal D_n\setminus T,\eta)^{\geo,p}\to \underset {S}\to {\varprojlim}\ \pi_1(X_{\overline K}\setminus (T\cup S)_{\overline K}, \eta)^{p}$, 
(resp. $\phi_n\defeq \underset{T} \to {\varprojlim}\ \phi_{n,T}: \underset{T}\to {\varprojlim}\ \pi_1(\Cal D_n\setminus T,\eta)^{\geo,p}\to 
\underset{(T,S)}\to {\varprojlim}\ \pi_1(X_{\overline K}\setminus (T\cup S)_{\overline K},\eta)^{p}$) makes  $\pi_1(\Cal D_n\setminus T,\eta)^{\geo,p}$
(resp. $\underset{T}\to {\varprojlim}\ \pi_1(\Cal D_n\setminus T,\eta)^{\geo,p}$) into a direct factor of
$\underset {S}\to {\varprojlim}\ \pi_1(X_{\overline K}\setminus (T\cup S)_{\overline K}, \eta)^{p}$ (resp. $\underset{(T,S)}\to {\varprojlim}\ 
\pi_1(X_{\overline K}\setminus (T\cup S)_{\overline K},\eta)^{p}$). In particular, both $\pi_1(\Cal D_n\setminus T,\eta)^{\geo,p}$ and 
$\underset{T}\to {\varprojlim}\ \pi_1(\Cal D_n\setminus T,\eta)^{\geo,p}$ are free pro-$p$ groups of infinite ranks.
Furthermore, the homomorphism $\psi_{n}:\pi_1(\Cal C_n,\eta)^{\geo,p}\to \underset {T}\to {\varprojlim}\ \pi_1(\Cal D_n\setminus T,\eta)^{\geo,p}$
makes $\pi_1(\Cal C_n,\eta)^{\geo,p}$ into a direct factor of $\underset {T}\to {\varprojlim}\ \pi_1(\Cal D_n\setminus T,\eta)^{\geo,p}$, and $\pi_1(\Cal C_n,\eta)^{\geo,p}$
is a free pro-$p$ group of infinite rank.
\endproclaim

\demo{Proof} First, we prove the assertion regarding the homomorphism $\phi_{n,T}$ by proving that 
 the map $H^1(\underset{S}\to {\varprojlim}\ \pi_1(X_{\overline K}\setminus (T\cup S)_{\overline K},\eta),\Bbb Z/p\Bbb Z)
\to H^1(\pi_1(\Cal D_n\setminus T,\eta)^{\geo},\Bbb Z/p\Bbb Z)$ induced by $\phi_{n,T}$ on cohomology is surjective, the result will then follow from 
Proposition 1.3.
We can (after possibly passing to a finite extension of $K$)  assume that $K$ contains the $p$-th roots of unity and all points in $T$ are $K$-rational. 
(In this, and other, proofs we will often use this argument. This is permissible because the various (pro-$p$, pro-$\ell$, and full) fundamental groups
under consideration are geometric fundamental groups.)
Let $f:\Cal U\to \Cal D_n$ be a generically $\mu_{p}$-torsor with $\Cal U$ geometrically connected and which is ramified above $T$.
Let $\{\Cal D_s^o\}_{s\in T}$ be pairwise disjoint open discs centred at the points $s\in T$, 
$\Cal X\defeq \Cal D_n\setminus (\cup_s \Cal D_s^o)$ an affinoid subdomain, and
$\tilde f:\Cal Y\to \Cal X$ the restriction of $f$ which is a $\mu_{p}$-torsor. When $\Cal D_s^o$
is small enough $f^{-1}(\Cal D_s^o)$ is geometrically connected as $f$ is totally ramified above $s$, $\forall s\in T$,  
and $\Cal Y$ is then geometrically connected, which we will assume form now on (cf. the argument in the proof of Proposition 2.3.2).
We can assume after possibly a finite extension of $R$ that $\Cal X$ has a (canonical) $R$-formal model
$Z=\Spf A$ with $Z_k$ reduced (cf. [Bosch-L\"utkebohmert-Raynaud], Theorem 1.3).

Let $h:Y\to Z$ be the finite morphism where $Y$ is the 
normalisation of $Z$ in $\Cal Y$. After possibly passing to a finite extension of $K$ we can assume that $Y_k$ is reduced (cf. [Epp]).
The $\mu_{p}$-torsor $\tilde f$ is given by a Kummer equation $\beta ^{p}=g$ where $g$ is a unit on $\Cal X$.
Let $V_K\defeq X_K\setminus (T\cup\{\infty\})$, then $(\Cal X,V_K)$ is a Runge pair.
The function $g$ can be approximated by a regular function $\tilde g$ on $V_K$ (cf. Proposition 2.2.2). For $\tilde g$ close to $g$ the equation 
$\beta^{p}=\tilde g$ defines a (possibly ramified) Galois covering $f_1:W_K\to X_K$ of degree $p$, with $W_K$ geometrically connected, 
whose pull-back to $\Cal X$ (via the rigid GAGA functor) is isomorphic to $\tilde f$. 
More precisely, one can write $g=\pi^tg_0$ where $g_0\in A$ is a unit and $0\le t<p$ an integer. One verifies easily that $t=0$ since $Y_k$ and $Z_k$ are reduced.
Let $\tilde g\in A^{\alg}$ such that $\tilde g-g\in \pi^r A$ where $U^{\alg}=\Spec A^{\alg}$ (cf. Propositions 2.2.1, Proposition 2.2.2, and the notations therein). 
Then for $r$ large enough $\tilde g g^{-1}\in 1+\pi^rA$ is a $p$-th power in $A$ and the Galois covering $f_1:Z_K\to X_K$ generically defined by the equation 
$\beta^{\ell}=\tilde g$ satisfies the above property. 
(More precisely, let $f\in A$, $r\ge 1$ large enough, then we can find $g\in A$ and $t\ge 1$ such that $1+\pi^rf=(1+\pi^tg)^{p}$.
If $r>\frac {pv_K(p)}{p-1}$, set $t=r-v_K(p)\ge 1$. Then $1+\pi^rf=1+p\pi^t(g+u_2\pi^tg^2+\ldots+u_{p}\pi^{t(p-1)-v_K(p)}g^{p})$, where $u_i\in R$ are units, and by Hensel's lemma we can find $g\in A$ such that $vf=g+u_2\pi^tg^2+\ldots+\pi^{t(p-1)}g^{p}$  where $v=\pi^{r-t}p^{-1}\in R$ is a unit (cf. [Bourbaki], Chapter III, $\S4.3$, Theorem 1).) 
Further, $f_1$ induces (via the rigid GAGA functor) 
a $\mu_{p}$-torsor $f_2:\Cal V\to \Cal X$ which is isomorphic to $\tilde f$, and a generically Galois cover $f_3:\Cal U\to X_K^{\rig}\setminus (\cup_s \Cal D_s^o)$.
One can then glue the generically $\mu_{p}$-torsors  $f$ and $f_3$ along $f_2$ and construct (via the rigid GAGA functor) a Galois covering 
$Y_K\to X_K$ which is ramified above $T$ and possibly a finite set 
$\Tilde S\subseteq X_K\setminus \Cal D_n$, and whose class in 
$H^1(\underset{S}\to {\varprojlim}\ \pi_1(X_{\overline K}\setminus (T\cup S)_{\overline K},\eta),\Bbb Z/\ell\Bbb Z)$ maps to the class of $f$ in
$H^1(\pi_1(\Cal D_n\setminus T,\eta)^{\geo},\Bbb Z/\ell\Bbb Z)$. The assertion regarding the homomorphism $\phi_n$ is proven in a similar way.
Similarly, using the fact that $(\Cal C_n,U_K)$ is a Runge pair, one proves that the natural map 
$H^1(\underset{(T,S)}\to {\varprojlim}\ 
\pi_1(X_{\overline K}\setminus (T\cup S)_{\overline K},\eta),\Bbb Z/p\Bbb Z)\to H^1(\pi_1(\Cal C_n,\eta)^{\geo},\Bbb Z/p\Bbb Z)$
is surjective, hence the second assertion follows from Proposition 1.3. 
Finally, the assertions on infinite rank follow from the fact that
$H^1(\pi_1(\Cal C_n,\eta)^{\geo},\Bbb Z/p\Bbb Z)$ is an infinite dimensional $\Bbb F_p$-vector space.
\qed
\enddemo

Recall the notations in Proposition 2.3.1, and 
consider the following exact sequence
$$1\to {\Cal H_n}\to  \pi_1(\Cal C_{n},\eta)^{\geo}  \to \Gamma \to 1,$$
where $\Cal H_n\defeq \Ker \lgroup \pi_1(\Cal C_{n},\eta)^{\geo}  \twoheadrightarrow \Gamma\rgroup$.
Further, let $P_n\defeq \Cal H_n^p$ 
be the maximal pro-$p$ quotient of
${\Cal H_n}$.
By pushing out the above sequence by the 
characteristic quotient  ${\Cal H_n}\twoheadrightarrow P_n$ we obtain an exact sequence
 $$1\to {P_n}\to \Delta_n \to \Gamma \to 1.$$
Similarly, consider the following exact sequence
$$1\to {\Cal H'_n}\to  \pi_1(\Cal D_{n}\setminus \{0\},\eta)^{\geo}  \to \widetilde \Gamma \to 1,$$
where $\Cal H'_n\defeq \Ker \lgroup \pi_1(\Cal D_{n}\setminus \{0\},\eta)^{\geo}  \twoheadrightarrow \widetilde \Gamma\rgroup$.
Further, let $\widetilde P'_n\defeq {(\Cal H'_n)}^p$ 
be the maximal pro-$p$ quotient of
$\Cal H'_n$.
By pushing out the above sequence by the 
characteristic quotient ${\Cal H'_n}\twoheadrightarrow \widetilde P'_n$ we obtain an exact sequence
 $$1\to \widetilde P'_n\to \widetilde \Delta '_n\to \widetilde \Gamma \to 1.$$

\proclaim {Proposition 2.3.4} Assume $K$ of equal characteristic $p\ge 0$. Then the natural morphism $\Cal C_n\to \Cal D_n$ induces a commutative diagram of exact sequences
$$
\CD
 1@>>> P_n @>>> \Delta_n @>>> \Gamma @>>> 1 \\
@.      @VVV     @VVV   @VVV \\
1@>>> \widetilde P'_n @>>> \widetilde \Delta'_n @>>> \widetilde \Gamma @>>> 1 \\
\endCD
$$
where the right vertical homomorphism $\Gamma\to \widetilde \Gamma$ is an isomorphism (cf. Lemma 2.3.1) 
and the middle vertical homomorphism $\Delta_n\to \widetilde \Delta'_n$ makes $\Delta_n$ into a semi-direct factor of $\widetilde \Delta'_n$ (cf. Lemma 1.5). 
\endproclaim

\demo{Proof} We only present the proof in the case $p>0$, the proof in the case $p=0$ is the same except for obvious simplifications. 
Let $\Delta_n \twoheadrightarrow G$ be a finite quotient which sits in an exact sequence $1\to Q\to G\to \Gamma_e\to 1$ where $\Gamma_e$
is the unique quotient of $\Gamma$ of cardinality $e$; for some integer $e$ 
prime-to-$p$ (cf. Proposition 2.3.1), with $Q$ a $p$-group. We will show there exists a surjective homomorphism
$\widetilde \Delta'_n\twoheadrightarrow G$ whose composition with $\Delta_n \to \widetilde \Delta'_n$ is the above homomorphism.
We can assume (without loss of generality) that the corresponding Galois covering $\Cal Y\to \Cal C_n$ with group $G$; 
$\Cal Y$ is normal and geometrically connected, is defined over $K$. This covering 
factorizes as $\Cal Y\to \Cal Y'\to \Cal C_{n}$ where  $\Cal Y'\to \Cal C_n$ is Galois with group 
$\Gamma _e\isom \mu_e$ and $\Cal Y\to \Cal Y'$ is Galois with group $Q$. 
After possibly a finite extension of $K$ we can assume that $n$ is divisible by $e$, the $\mu_e$-torsor $\Cal Y'\to \Cal C_n$
is generically defined by an equation $\widetilde Z^e=Z$ for a suitable choice of the parameter $Z$, and 
$\Cal Y'=\Cal C_{\frac {n}{e}}=\Sp \frac {K<\widetilde Z,\widetilde Y,\widetilde W>}{(\widetilde Z-\pi^{\frac {n}{e}}\widetilde Y,
\widetilde Y\widetilde W-1)}$ ($\widetilde Y^e=Y$) is an annulus of thickeness $0$. 
The $\mu_e$-torsor $\Cal Y'\to \Cal C_n$ extends then to a generically $\mu_e$-torsor $\Cal X'\to \Cal D_n$ generically defined by an equation $\widetilde Z^e=Z$, which 
is (totally) ramified only above $0$, and $\Cal X'=\Cal D_{\frac {n}{e}}=
\Sp \frac {K<\widetilde Z,\widetilde Y>}{(\widetilde Z-\pi^{\frac {n}{e}}\widetilde Y)}$ is a closed disc centred at the unique point above $0\in \Cal D_n$; which we denote also $0$.

By Proposition 2.3.2 applied to $\Cal Y'\to \Cal X'$ there exists (after possibly a finite extension of $K$) a Galois 
covering $\Cal X\to \Cal X'$ with group $Q$, ramified only above $0$, 
with $\Cal X$ geometrically connected and such that we have a commutative diagram of cartesian squares.
$$
\CD
\Cal Y @>>> \Cal Y' @>>> \Cal C_n\\
@VVV @VVV   @VVV\\
\Cal X @>>> \Cal X' @>>> \Cal D_n \\
\endCD
$$

Next, we borrow some ideas from [Garuti] (preuve du Th\'eor\`eme 2.13). We claim one can choose the above (geometric) covering $\Cal X\to \Cal X'$ 
such that the composite covering $\Cal X\to \Cal D_n$ is Galois with group $G$. Indeed, consider the quotient
$\Delta_n \twoheadrightarrow \Delta _{\Cal Y'}$ (resp. $\widetilde \Delta_n' \twoheadrightarrow \Delta _{\Cal X'}$)
of $\Delta_n$ (resp. $\widetilde \Delta_n'$) which sits in the following exact sequence $1\to P_{\Cal Y'}\to \Delta_{\Cal Y'}\to \Gamma_e \to 1$
(resp. $1\to P_{\Cal X'}\to \Delta_{\Cal X'}\to \widetilde \Gamma_e \to 1$) where $P_{\Cal Y'}\defeq \pi_1(\Cal Y',\eta)^{\geo,p}$ (resp. $P_{\Cal X'}\defeq
\pi_1(\Cal X'\setminus \{0\},\eta)^{\geo,p}$) and $\widetilde \Gamma_e$ is the unique quotient of $\widetilde \Gamma$ of cardinality $e$. 
We have a commutative diagram of exact sequences

$$
\CD
1 @>>> P_{\Cal Y'} @>>> \Delta_{\Cal Y'} @>>> \Gamma_e @>>>1\\
@. @VVV @VVV @VVV\\
1 @>>> P_{\Cal X'} @>>> \Delta_{\Cal X'} @>>> \widetilde \Gamma_e @>>>1\\
\endCD
$$
where the right vertical map is an isomorphism (cf. Proposition 2.3.1). 
The choice of a splitting of the upper sequence in the above diagram (which splits since $P_{\Cal Y'}$ is pro-$p$ and $\Gamma _e$ is (pro-)prime-to-$p$) induces
an action of $\Gamma_e$ on $P_{\Cal X'}$, and $P_{\Cal Y'}$ is a direct factor of $P_{\Cal X'}$ (cf. Proposition 2.3.2) which is stable by this action of $\Gamma_e$. Further, $P_{\Cal Y'}$ has a supplement $E$ in $P_{\Cal X'}$ which is invariant under the action of $\Gamma_e$ by [Garuti], Corollaire 1.11. The existence of this supplement $E$ implies that one can choose 
$\Cal X\to \Cal X'$ as above such that the finite composite covering $\Cal X\to \Cal D_n$ is Galois with group $G$.
More precisely, if the Galois covering $\Cal Y\to \Cal Y'$ corresponds to the surjective homomorphism $\rho:P_{\Cal Y'}\twoheadrightarrow Q$ 
(which is stable by $\Gamma_e$ since $\Cal Y\to \Cal C_n$ is Galois) then we consider the Galois covering   $\Cal X\to \Cal X'$ 
corresponding to the surjective homomorphism $\widetilde P_{\Cal X'}=P_{\Cal Y'} \star E\twoheadrightarrow  Q$ 
which is induced by $\rho$ and the trivial homomorphism $E\to Q$, which is stable by $\Gamma_e$. 

The above construction can be performed in a functorial way with respect to the various finite quotients of $\Delta_n$. More precisely,
let $\{\phi_j:\Delta_n \twoheadrightarrow G_j\}_{j\in J}$ be a cofinal system of finite quotients of $\Delta_n$ where $G_j$
sits in an exact sequence $1\to Q_j\to G_j\to \Gamma_{e_j}\to 1$, for some integer $e_j$ 
prime-to-$p$, and $Q_j$ a $p$-group. Assume we have a factorisation $\Delta_n \twoheadrightarrow G_{j'}\twoheadrightarrow G_j$ for $j',j\in J$. Thus, $e_j$ divides $e_{j'}$, and we can assume without loss of generality (after replacing the group extension $G_j$ by its pull-back via $\Gamma _{e_j'}\twoheadrightarrow \Gamma _{e_j}$)
that $e\defeq e_j=e_{j'}$. With the above notations we then have surjective homomorphisms $\rho_{j'}:P_{\Cal Y'}\twoheadrightarrow Q_{j'}$, 
$\rho_j:P_{\Cal Y'}\twoheadrightarrow Q_j$ 
(which are stable by $\Gamma_e$), and $\rho_{j}$ factorises through $\rho_{j'}$. Then we consider the Galois covering(s)   $\Cal X_{j'}\to \Cal X'$ (resp. 
$\Cal X_{j}\to \Cal X'$) corresponding to the surjective homomorphism(s) $\psi_{j'}:P_{\Cal X'}=P_{\Cal Y'} \star E\twoheadrightarrow  Q_{j'}$ 
(resp. $\psi_j:P_{\Cal X'}=P_{\Cal Y'} \star E\twoheadrightarrow  Q_j$) 
which are induced by $\rho_{j'}$ (resp. $\rho_j$) and the trivial homomorphism $E\to Q$, which are stable by $\Gamma_e$ and $\psi_{j}$ factorises through
$\psi_{j'}$. Finally, we deduce from this construction the existence of a surjective continuous homomorphism $\widetilde \Delta_n'\twoheadrightarrow  \Delta_n$ which is a left 
inverse to the natural homomorphism $\Delta_n \to \widetilde \Delta_n'$ (cf. Lemma 1.5).
\qed
\enddemo

Next, recall the discussion and notations before Proposition 2.3.3, and consider the following exact sequence 
$$1\to \widetilde {\Cal H}_n\to  \underset{T} \to {\varprojlim}\ \pi_1(D_{n}\setminus T,\eta)^{\geo}  \to \widetilde \Gamma \to 1,$$
where $\widetilde {\Cal H}_n\defeq \Ker \lgroup \underset{T} \to {\varprojlim}\ \pi_1(D_{n}\setminus T,\eta)^{\geo}  \twoheadrightarrow \widetilde \Gamma \rgroup$.
Further, let $\widetilde P_n\defeq {\widetilde {\Cal H}_n}^p$ 
be the maximal pro-$p$ quotient of
$\widetilde {\Cal H}_n$.
By pushing out the above sequence by the 
characteristic quotient $\widetilde {\Cal H}_n\twoheadrightarrow \widetilde P_n$ we obtain an exact sequence
 $$1\to \widetilde {P}_n\to \widetilde \Delta_n\to \widetilde \Gamma \to 1.$$

\proclaim {Proposition 2.3.5} Assume $\char (K)=0$ with no restriction on $\char(k)=p\ge 0$. Then the morphism $\Cal C_n\to \Cal D_n$ induces a commutative diagram of exact sequences
$$
\CD
 1@>>> P_n @>>> \Delta_n @>>> \Gamma @>>> 1 \\
@.      @VVV     @VVV   @VVV \\
1@>>> \widetilde {P}_n @>>> \widetilde \Delta_n @>>> \widetilde \Gamma @>>> 1 \\
\endCD
$$
where the right vertical homomorphism $\Gamma\to \widetilde \Gamma$ is an isomorphism (cf. Lemma 2.3.1) 
and the middle vertical homomorphism $\Delta_n\to \widetilde \Delta_n$ makes $\Delta_n$ into a semi-direct factor of $\widetilde \Delta_n$ (cf. Lemma 1.5). 
\endproclaim

\demo{Proof} We only explain the proof in the case $p>0$, the proof in the case $p=0$ is the same except for obvious simplifications. 
The proof is entirely similar to the proof of Proposition 2.3.4 using Proposition 2.3.3 instead of Proposition 2.3.2. 
With the notations in the proof of loc. cit. one applies Proposition 2.3.3 to $\Cal Y'\to \Cal X'$
to ensure the existence (after possibly a finite extension of $K$) of a Galois 
covering $\Cal X\to \Cal X'$ with group $Q$, ramified above a finite set $T'\subset \Cal X'$, 
with $\Cal X$ geometrically connected and such that we have a commutative diagram as in loc. cit.
One then considers the quotients $\Delta_n\twoheadrightarrow \Delta_{\Cal Y'}$ as in loc. cit., 
and $\widetilde \Delta_n \twoheadrightarrow \Delta _{\Cal X'}$
which sits in the following exact sequence $1\to P_{\Cal X'}\to \Delta_{\Cal X'}\to \widetilde \Gamma_e \to 1$ where $P_{\Cal X'}\defeq
\pi_1(\Cal X'\setminus T',\eta)^{\geo,p}$ and follow the same arguments as in loc. cit.
\qed
\enddemo

\definition{Remark 2.3.6} With the same notations as in Propositions 2.3.4 and 2.3.5 the pro-$p$ group $\widetilde P_n'$ (resp. $\widetilde P_n$) is free
and the homomorphism  $P_n\to \widetilde P_n'$ (resp. $P_n\to \widetilde P_n$) makes $P_n$ into a direct factor of $\widetilde P_n'$ (resp. $\widetilde P_n$).
\enddefinition

\definition{Remark 2.3.7} The proofs of Propositions 2.3.4 and 2.3.5 rely on Corollaire 1.11 in [Garuti] which relies on Lemme 1.9 and Lemme 1.10 in loc. cit.  
We take this opportunity to make precise some steps in the proof of these Lemmas. Using the notations in loc. cit., in the proof of Lemme 1.9; the definition of the set of pairs
$(S,\varphi_S)$, one should require $S\subseteq L^{\star}$ and the given action $\varphi_S$ of $\Gamma$ on $L/S$ to lift the given action of $\Gamma$ on $L/L^{\star}$.
Also the group $\Gamma_1$ should be defined as the subgroup of $\Aut(L)$ generated by lifts of the elements of $\Gamma$ viewed as acting on $L/N$ (automorphisms of 
$L/N$ which lift automorphisms of $L/L^{\star}$ lift to automorphisms of $L$ cf. Corollaire 1.7 in loc. cit.) 
so that $\Gamma_1$ acts on $L/(N\cap U)$ stabilising $N/(N\cap U)$,
$\Gamma_2$ should then be defined by first taking the image of $\Gamma_1$ in $\Aut (L/(N\cap U))$ (this image maps to $\Aut (L/N)$) and then take the pull back of this 
image via the map $\Gamma\twoheadrightarrow \Im (\Gamma_1\to \Aut(L/N))$.
Also in the proof of Lemme 1.10, the definition of the pair $(S,\alpha_S)$, with the notations in loc. cit. one should require $S\subseteq L^{\star}$,
%the given actions $\Phi _S$ and $\Psi _S$ of $\Gamma$ on $L/S$ should lift the given action of $\Gamma$ on $L/L^*$,
and the automorphism $\alpha_S\in \Aut (L/S)$ should lift the identity of $\Aut (L/L^{\star})$. 
\enddefinition

\subhead
\S3 Geometric fundamental groups of affinoid $p$-adic curves
\endsubhead
In this section we investigate the structure of the geometric fundamental group of rigid affinoid $K$-curves which are embedded in a proper $K$-curve.
Let $X$ be a proper and normal formal $R$-curve with $X_K$ smooth,
$U\hookrightarrow X$ an $R$-formal affine sub-scheme, and $\Cal U\defeq U_K\hookrightarrow X_K$ the associated $K$-rigid analytic affinoid space 
(which is an affinoid rigid subspace of $X_K$). 
We assume that the special fibre $U_k$ 
of $U$ is connected, reduced, and $X_k\setminus U_k=\{\bar x_i\}_{i=1}^m$ consists of a finite set of closed points 
where $\bar x_i\in X_k(k)$ is a smooth point of $X_k$, $1\le i\le m$. 
Let ${\Cal F}_{i}\defeq \Spf \hat \Cal O_{X,{\bar x_i}}$ be the formal germ of $X$ at the point $\bar x_i$, $1\le i\le m$. Thus, 
$\hat \Cal O_{X,{\bar x_i}}\isom R[[T_i]]$ is the completion of the local ring $\Cal O_{X,{\bar x_i}}$. 
Let $\Cal D_i\defeq \Sp K<T_i>$ be the rigid closed unit disc and $\Cal C_i\defeq \Sp K<T_i,\frac {1}{T_i}>$
the rigid standard annulus of thickness $0$. The formal fibre of $\bar x_i$ in $X_K$ is isomorphic to the open disc $\Cal D_i^o\defeq \Cal D_i\setminus \Cal C_i$
and $X_K\setminus U_K$ is isomorphic to the disjoint union $\bigcup_{i=1}^m\Cal D_i^o$.
In what follows we identify the formal fibre of $\bar x_i$ in $X_K$ with the open disc $\Cal D_i^o$, $1\le i\le m$.
Write $x_i\in \Cal D_{i,K}(K)$ for the zero point $T_i=0$ of $\Cal D_{i,K}$ which we view, via the above identification, as a point in $X_K(K)$. 
The above conditions on the affinoid curve $\Cal U$ are not too restrictive. More precisely, we have the following.

\proclaim {Theorem 3.1} Let $Y_K$ be a smooth and geometrically connected affinoid $K$-curve. Then, after possibly a finite extension of $K$, one can
embed $Y_K$ into a proper, geometrically connected, and smooth $K$-curve $X_K$ such that the complement $X_K\setminus Y_K$ consists of a disjoint union of finitely many 
open unit $K$-discs as in the above discussion where $Y_K=\Cal U$.
\endproclaim
\demo{Proof} See [Van Der Put], Theorem 1.1.
\qed
\enddemo

We use the notations in $2.3$. For an integer $n\ge0$ we write $\Cal V_n\defeq X_K\setminus \lgroup \bigcup_{i=1}^m \{\Cal D_{i,n}^o\}\rgroup$ where $\Cal D_{i,n}\defeq \Cal D_n\subseteq \Cal D_i$ and $\Cal D_{i,n}^o\defeq \Cal D_n^o\defeq \Cal D_n\setminus \Cal C_n$ are as in loc. cit. 
Thus, $\Cal U\subseteq \Cal V_n$ is a quasi-compact rigid analytic subspace of $X_K$. We also write $\Cal C_{i,n}\defeq \Cal C_n\subset \Cal D_{i,n}$ which is a closed annulus of thickness $0$.

\proclaim{Proposition 3.2} Let $S\subset \Cal U$ be a finite set of points and $f:\Cal Z\to \Cal U$ a finite Galois covering with Galois group $G$ which is \'etale above $\Cal U\setminus S$.
Then, after possibly a finite extension of $K$, there exists $n>0$ such that $f$ extends to a finite Galois covering $f_n:\Cal Z_n\to \Cal V_n$ which is Galois with Galois group $G$
and \'etale above $\Cal U\setminus S$. Moreover, let $f_i\defeq f_{i,n}:\Cal W_i=\bigcup_j\Cal W_{i,j}\to \Cal C_{i,n}$ be the restriction of $f_n$ to the annulus $\Cal C_{i,n}$, 
where $\{\Cal W_{i,j}\}_j$ are the connected components of $\Cal W_i$, $1\le i\le m$. Then the decomposition group $G_{i,j}\subseteq G$ of each connected component 
$\Cal W_{i,j}$ is a solvable group which is an extension of a cyclic group of order prime-to-$p$ by a $p$-group.
\endproclaim

\demo {Proof} This can be proven using similar arguments used by Raynaud in [Raynaud] to prove a similar result in the case where $\Cal U$ is the closed unit disc centred at $0$
which is embedded in $(\Bbb P^1_K)^{\rig}$ as the complement of the open disc centred at $\infty$ (see Remarques 3.4.12(i) in loc. cit.). 
We briefly explain the outline of proof. First, one can assume (without loss of generality) that 
$S=\emptyset$, $m=1$, $x\defeq x_1$, and $\Cal D\defeq \Cal D_1$. 
%and replace $\Cal U$ by a Laurent subdomain $\Cal U(\frac {1}{f})\defeq \{z\in \Cal U\ :\ \vert f(z)\vert\ge 1\}$
%where $f$ is a local coordinate at $x$. 
Using similar arguments as in [Raynaud] Proposition 3.4.1 one can extend the \'etale covering $f$ to an \'etale covering $f_{n'}:\Cal Z_{n'}\to \Cal U_{n'}\defeq
X_K\setminus \Cal D^o_{\frac {1}{n'}}$ where $\Cal D^o_{\frac{1}{n'}}\defeq \{x\in \Cal D\ :\ \vert Z(x) \vert< \vert \pi\vert ^{\frac {1}{n'}}\}$ for some integer $n'\ge 1$ and the germ of such an extension is unique. (More precisely, one establishes the analogs of 
Lemme 3.4.2 and Lemme 3.4.3 in loc. cit. using similar arguments.) In particular, there exists $n'$ as above such that $f_{n'}$ is Galois with Galois group $G$ (cf. loc. cit. Corollaire 3.4.8).
There exists a formal model $U_{n'}$ of $\Cal U_{n'}$ whose (reduced) special fibre consists of the special fibre $U_k$ of $\Cal U$ which is linked to an affine line at a double point $\bar x$ in which specialise the points of the open annulus $\Cal A_{\frac{1}{n'}}\defeq \{x\in \Cal D\ :\ 1>\vert Z(x)\vert>\vert \pi\vert ^{\frac {1}{n'}}\}$ 
(cf. [Raynaud] 3.3.3 in the special case where $\Cal U$ is the unit closed disc). 
Let $g_{n'}:Z_{n'}\to U_{n'}$ be the morphism of normalisation of $U_{n'}$ in 
$\Cal Z_{n'}$ which is Galois with group $G$, $g_{n'}^{-1}(U_k)$ the reduced inverse image of $U_k$ in $Z_{n'}$, and $\{y_i\}_i$ the points of $g_{n'}^{-1}(U_k)$ above $\bar x$. 
Then one proves that for $n'$ large enough $g_{n'}^{-1}(U_k)$ is normal at the points $\{y_i\}_i$ (cf. loc. cit. Lemme 3.4.2), 
there is a one-to-one correspondence between the $\{y_i\}_i$ and the connected components  $\{\Cal Y_i\}_i$ of the inverse image $f_{n'}^{-1}(\Cal A_{\frac{1}{n'}})$ 
in $\Cal Z_{n'}$ of the open annulus  
$\Cal A_{\frac {1}{n'}}$, as well as a one-to-one correspondence between the corresponding decomposition groups (cf. loc. cit. Proposition 3.4.6(ii)). Moreover, the decomposition group of such a component $\Cal Y_i$ is solvable (cf. loc. cit. Corollaire 3.4.8).
Let $Z$ be the normalisation of $U$ in $\Cal Z$. After possibly a finite extension of $K$ we can assume that the special fibre $Z_k$ of 
$\Cal Z$ is reduced (cf. [Epp]). In this case the decomposition group of a connected component $\Cal Y_i$ as above is an extension of a cyclic group of order prime-to-$p$ by a $p$-group.
Indeed, in this case with the notations of loc. cit., the proof of Corollaire 3.4.8, the group $I_i$ is a $p$-group as we assumed $Z_k$ is reduced. Finally, after a finite extension of $K$ 
we can assume the above open annulus $\Cal A_{\frac{1}{n'}}$ of thickness $\frac {1}{n'}$ (a rational) is an annulus $\{x\in \Cal D\ :\ 1>\vert Z(x)\vert>\vert \pi\vert ^{n}\}$ of thickness $n$ for some integer $n>0$.
\qed
\enddemo

For the remaining of this section, let $S\subset \Cal U$ be a (possibly empty) {\it finite} set of points.

\proclaim {Theorem 3.3} Assume $K$ of equal characteristic $p>0$, and let $\ell$ be a prime integer (possibly equal to $p$). 
Then the morphism $\Cal U\to X_K$ induces (via the rigid GAGA functor) a continuous homomorphism 
$\pi_1(\Cal U\setminus S,\eta)^{\geo}
\to \pi_1(X_K\setminus (\{x_i\}_{i=1}^m\cup S),\eta)^{\geo}$ (resp. $\pi_1(\Cal U\setminus S,\eta)^{\geo,\ell}
\to \pi_1(X_K\setminus (\{x_i\}_{i=1}^m\cup S),\eta)^{\geo,\ell}$)
which makes $\pi_1(\Cal U\setminus S,\eta)^{\geo}$ (resp. $\pi_1(\Cal U\setminus S,\eta)^{\geo,\ell}$) into a semi-direct factor of 
$\pi_1(X_K\setminus (\{x_i\}_{i=1}^m\cup S),\eta)^{\geo}$ (resp. direct factor of $\pi_1(X_K\setminus (\{x_i\}_{i=1}^m\cup S),\eta)^{\geo,\ell}$).
Moreover,  $\pi_1(\Cal U\setminus S,\eta)^{\geo,\ell}$ is a free pro-$\ell$ group of infinite (resp. finite) rank if $\ell=p$ (resp. if $\ell\neq p$). 
\endproclaim

\demo{Proof} We show the criterion in Lemma 1.5 is satisfied. 
Let $\pi_1(\Cal U\setminus S,\eta)^{\geo}\twoheadrightarrow G$ be a finite quotient which we assume (without loss of generality) corresponding to 
a finite Galois covering $f:\Cal V\to \Cal U$ with group G, 
\'etale above $\Cal U\setminus S$, with $\Cal V$ normal and geometrically connected. 
We will show the existence of a surjective homomorphism $\pi_1(X_K\setminus (\{x_i\}_{i=1}^m\cup S),\eta)^{\geo}\twoheadrightarrow G$ whose composite with  $\pi_1(\Cal U\setminus S,\eta)^{\geo}
\to \pi_1(X_K\setminus (\{x_i\}_{i=1}^m\cup S),\eta)^{\geo}$ is the above homomorphism. We assume the existence of an extension $f_n:\Cal Z_n\to \Cal V_n$ of $f$ as in Proposition 3.2. 
For $1\le i\le m$, let $f_i\defeq f_{i,n}:\Cal W_i=\bigcup_{j=1}^{t_i}\Cal W_{i,j}\to \Cal C_{i,n}$ be the restriction of $f_n$  to the annulus $\Cal C_{i,n}$ with $\{\Cal W_{i,j}\}_{j=1}^{t_i}$ the connected components of $\Cal W_i$, and $G_{i,j}\subseteq G$ the decomposition group of 
$\Cal W_{i,j}$ which is an extension of a cyclic group of order prime-to-$p$ by a $p$ group.
Fix $1\le j_0\le t_i$, then $f_i\isom \Ind _{G_{i,j_0}}^G f_{i,j_0}$ is an induced covering (cf. [Raynaud 4.1) where $f_{i,j_0}:\Cal W_{i,j_0}\to \Cal C_{i,n}$ is the restriction of $f_i$
to $\Cal W_{i,j_0}$.

By Proposition 2.3.4 (the equal characteristic $p>0$ case) there exists (after possibly a finite extension of $K$) a finite Galois covering $\tilde f_{i,j_0}:\Cal Y_{i,j_0}\to \Cal D_{i,n}$ with group $G_{i,j_0}$, $\Cal Y_{i,j_0}$ is normal and geometrically 
connected, whose pull-back to ${\Cal C}_{i,n}$ via the embedding ${\Cal C}_{i,n}\hookrightarrow \Cal D_{i,n}$ is isomorphic to $f_{i,j_0}$, and 
$\tilde f_{i,j_0}$ is ramified only above $x_i$.
Let $\tilde f_i : \Cal Y_i\defeq \Ind_{G_{i,j_0}}^G \Cal Y_{i,j_0}\to \Cal D_{i,n}$ be the induced coverings (cf. loc. cit.), $1\le i\le m$. 
One can then patch the covering $f_n$ with the coverings $\{\tilde f_i\}_{i=1}^m$ along the restrictions of these coverings above the annuli $\Cal C_{i,n}$ 
(the restriction of $f_n$ and $\tilde f_{i,n}$ to $\Cal C_{i,n}$ are isomorphic by construction)
to construct a finite Galois covering  $\tilde f:Y_K\to X_K$ between smooth and proper rigid $K$-curves with group $G$, $Y_K$ is geometrically connected, 
which gives rise (via the rigid GAGA functor) to a homomorphism 
$\pi_1(X_K\setminus (\{x_i\}_{i=1}^n\cup S),\eta)^{\geo}\twoheadrightarrow G$ as required.

Moreover, one verifies easily that the above construction can be performed in a functorial way with respect to the various quotients 
of $\pi_1(\Cal U\setminus S,\eta)^{\geo}$ (in the sense of Lemma 1.5) using Proposition 2.3.4, so that it induces a continuous homomorphism
$\pi_1(X_K\setminus (\{x_i\}_{i=1}^m\cup S),\eta)^{\geo}\to \pi_1(\Cal U\setminus S,\eta)^{\geo}$ which is left inverse to $\pi_1(\Cal U\setminus S,\eta)^{\geo}
\to \pi_1(X_K\setminus (\{x_i\}_{i=1}^m\cup S),\eta)^{\geo}$. The second assertion is proven in a similar way.
Note that $\pi_1(X_K\setminus (\{x_i\}_{i=1}^m\cup S),\eta)^{\geo,\ell}$ is pro-$\ell$ free (cf. [Serre1], Proposition 1, and Proposition 1.1.1, in the case $\ell=p$, 
and [Grothendieck], Expos\'e XIII, Corollaire 2.12, otherwise), the assertion that $\pi_1(\Cal U\setminus S,\eta)^{\geo,\ell}$ is free follows then from the above
discussion. Finally, the assertion on the rank follows from Proposition 3.5 below if $\ell=p$, and from the fact that $\pi_1(X_K\setminus (\{x_i\}_{i=1}^m\cup S),\eta)^{\geo,\ell}$
is finitely generated if $\ell\neq p$ (cf. loc. cit.)
\qed
\enddemo

%\proclaim {Proposition 3.4} Assume $\char (K)=p>0$. Then $\pi_1(\Cal U\setminus S,\eta)^{\geo,p}$ is a free pro-$p$ group of infinite rank. 
%\endproclaim

%\demo{Proof} The first assertion follows from Theorem 3.3. Let $\Cal U=\Sp \Cal A$. For the second assertion it suffices to show that 
%$H_{\et}^1(\Spec \Cal A_{\overline K},\Bbb Z/p\Bbb Z)$ is infinite ($\Cal A_{\overline K}\defeq \Cal A\otimes _K\overline K$) which follows easily
%from the structure of $\Cal A$ as an affinoid algebra.
%\qed
%\enddemo

Let $T\subset \bigcup_{i=1}^m \Cal D_{i}^o$ be a {\it finite} set of closed points of 
$X_K$. We view $T\subset X_K$ as a closed subscheme of $X_K$. We have an 
exact sequence of profinite groups
$$1\to \pi_1(X_{K}\setminus (T\cup S),\eta)^{\geo}\to \pi_1(X_K\setminus (T\cup S),\eta)\to 
\Gal (\overline K/K)\to 1.$$
By passing to the projective limit over all finite sets of closed points 
$T\subset  \bigcup_{i=1}^m \Cal D_{i}^o$ we obtain an exact sequence
$$1\to \underset{T} \to {\varprojlim}\ \pi_1(X_{K}\setminus (T\cup S),\eta)^{\geo}
\to \underset {T} \to {\varprojlim}\ \pi_1(X_{K}\setminus (T\cup S),\eta)\to \Gal (\overline 
K/K)\to 1.$$
The profinite group $\underset{T} \to {\varprojlim}\ \pi_1(X_{K}\setminus (T\cup S),\eta)^{\geo}$
is {\it free} if $\char (K)=0$ as follows from the well-known structure of the geometric \'etale fundamental groups of (affine) curves
in characteristic zero (cf. [Grothendieck], Expos\'e XIII, Corollaire 2.12).

\proclaim {Theorem 3.4}
Assume $\char (K)=0$ with no restriction on $\char(k)=p\ge 0$.
Let $\ell$ be a prime integer (possibly equal to $\char(k)$ if $\char(k)>0$). Then the morphism $\Cal U\to X_K$ induces 
(via the rigid GAGA functor) a continuous 
homomorphism $\pi_1(\Cal U\setminus S,\eta)^{\geo}
\to \underset{T} \to {\varprojlim}\  \pi_1(X_K\setminus (T\cup S),\eta)^{\geo}$ 
(resp. $\pi_1(\Cal U\setminus S,\eta)^{\geo,\ell}
\to \underset{T} \to {\varprojlim}\  \pi_1(X_K\setminus (T\cup S),\eta)^{\geo,\ell}$) 
which makes $\pi_1(\Cal U\setminus S,\eta)^{\geo}$ (resp. $\pi_1(\Cal U\setminus S,\eta)^{\geo,\ell}$) 
into a semi-direct factor of $\underset{T} \to {\varprojlim}\ \pi_1(X_K\setminus (T\cup S),\eta)^{\geo}$
(resp. direct factor of $\underset{T} \to {\varprojlim}\ \pi_1(X_K\setminus (T\cup S),\eta)^{\geo,\ell}$).
In particular, the pro-$\ell$ group $\pi_1(\Cal U\setminus S,\eta)^{\geo,\ell}$ is free.
\endproclaim

\demo{Proof} The proof is similar, almost word by word, to the proof of Theorem 3.3. One has to use Proposition 2.3.5  instead of the use of Proposition 2.3.4 made in the 
proof of Theorem 3.3.
\qed
\enddemo

\proclaim {Proposition 3.5} Assume $\char k=p>0$ with no restriction on $\char(K)$. Then the pro-$p$ group $\pi_1(\Cal U\setminus S,\eta)^{\geo,p}$ is free 
of infinite rank.
\endproclaim

\demo{Proof} The first assertion follows from Theorem 3.3 (resp. Theorem 3.4) if $\char(K)=p$ (resp. $\char(K)=0$). 
For the second assertion it suffices to show that the $\Bbb F_p$-vector space $H_{\et}^1(\Spec \Cal A_{\overline K},\Bbb Z/p\Bbb Z)$, where
$\Cal U=\Sp \Cal A$ and $\Cal A_{\overline K}\defeq \Cal A\otimes _K\overline K$, is infinite, which follows easily from the structure of $\Cal A$ as an affinoid algebra. 
\qed
\enddemo

\proclaim {Proposition 3.6} Let $\char(k)=p\ge 0$ with no restrictions on $\char(K)$.
Then the morphism $\Cal U\to X_K$ induces a continuous homomorphism $\pi_1(\Cal U\setminus S,\eta)^{\geo,p'}
\to \pi_1(X_K\setminus (\{x_i\}_{i=1}^m\cup S),\eta)^{\geo,p'}$ which makes $\pi_1(\Cal U\setminus S,\eta)^{\geo,p'}$ into a semi-direct factor of
$\pi_1(X_K\setminus (\{x_i\}_{i=1}^m\cup S),\eta)^{\geo,p'}$.
\endproclaim

\demo {Proof} The proof follows by using similar arguments to the ones used in the proofs of Theorems 3.3 and 3.4. 
More precisely, with the notations in the proofs of loc. cit. the morphism $\Cal W_{i,j}\to {\Cal C}_{i,n}$ in this case is a $\mu_e$-torsor, where 
$e$ is an integer prime-to-$p$, and extends (uniquely, after possibly a finite extension of $K$) to a cyclic Galois cover $\Cal Y_{i,j}\to \Cal D_{i,n}$ of 
degree $e$ ramified only above $x _i$ (cf. Lemma 2.3.1 and the isomorphism $\Gamma \isom \widetilde \Gamma$ therein). 
\enddemo

In what follows let $g\defeq g_{X_K}$ be the arithmetic genus of $X_K$ ($g_{\Cal U}\defeq g$ 
is also called the genus of the affinoid $\Cal U$).

\proclaim {Theorem 3.7} Let $\char(k)=p\ge 0$ with no restriction on $\char(K)$.
Let $S(\overline K)=\{y_1,\cdots,y_r\}$ of cardinality $r\ge 0$.
Then the homomorphism $\pi_1(\Cal U\setminus S,\eta)^{\geo,p'}
\to \pi_1(X_K\setminus (\{x_i\}_{i=1}^m\cup S),\eta)^{\geo,p'}$ (as in Proposition 3.6) is an isomorphism. In particular,
$\pi_1(\Cal U\setminus S,\eta)^{\geo,p'}$ is (pro-)prime-to-$p$ free on $2g+m+r-1$ generators and further can be generated by $2g+m+r$ generators 
$\{a_1,\cdots,a_g,b_1,\cdots,b_g,\sigma_1,\cdots,\sigma_m,\tau_1,\cdots,\tau_r\}$
subject to the unique relation $\prod_{j=1}^g[a_j,b_j]\prod _{i=1}^m\sigma _i\prod_{t=1}^r \tau _t=1$, where $\sigma_i$ (resp $\tau _t$) is a generator of inertia at $x _i$
(resp. $y_t$). 
\endproclaim

\demo{Proof} The second assertion follows from [Grothendieck], Expos\'e XIII, Corollaire 2.12. 
The homomorphism $\pi_1(\Cal U\setminus S,\eta)^{\geo,p'}\to \pi_1(X_K\setminus (\{x_i\}_{i=1}^n\cup S),\eta)^{\geo,p'}$ is injective by Proposition 3.6.
We show it is surjective. To this end it suffices to show that given a finite Galois covering $f:Y\to X$ with group $G$ of cardinality prime-to-$p$, with $Y$ normal and 
geometrically connected, which is \'etale above $X_K\setminus (\{x_i\}_{i=1}^m\cup S)$, and $\tilde f:\Cal V\to \Cal U$ its restriction to $\Cal U$, then $\Cal V$ is 
geometrically connected. 
We can assume, without loss of generality, that $Y_k$ is reduced (cf. Abhyankar's Lemma, [Grothendieck], Expos\'e X, Lemme 3.6).
First, note that $f^{-1}(\Cal D_i^o)$ is a disjoint union of finitely many formal open unit discs (cf. [Raynaud, Lemma 6.3.2), $1\le i\le m$. 
Let $V$ be the normalisation of $U$ in $\Cal V$.
Suppose that $\Cal V$ is disconnected, then $V_k$ is disconnected, and a fortiori $Y_k$ is also disconnected as $Y_k\setminus V_k$ is regular (cf. loc. cit.), 
but this contradicts the fact that $Y_K$ is connected. 
\qed
\enddemo

\definition {Remark 3.8} {\bf (i)}\ If $\char (k)=0$ the profinite group $\pi_1(\Cal U\setminus S,\eta)^{\geo}$ is {\it free} and finitely generated as follows from Theorem 3.7.
Apart from this case the profinite group $\pi_1(\Cal U\setminus S,\eta)^{\geo}$ is {\it not free} (neither is it finitely generated) as the ranks of its maximal pro-$\ell$ quotients can be different for different primes $\ell$ (cf. Theorems 3.3, 3.4, and 3.7).  In this sense Theorem 3.3 and Theorem 3.4 are optimal results one can prove regarding the structure of the 
{\it full} geometric fundamental group of a $p$-adic affinoid curve. 

{\bf (ii)}\ There is no analog in mixed characteristics to Theorem 3.4, for the full $\pi_1^{\geo}$, where one replaces the infinite union of the finite sets of points $T$ (as in loc. cit.) by a single fixed finite set of points $\widetilde T\subset \bigcup_{i=1}^m \Cal D_{i}^o$. More precisely, in this case one can not control the ramification arising from an extension 
of an \'etale covering $\Cal V\to \Cal U$ of the affinoid $\Cal U$ to a ramified covering $Y_K\to X_K$. 
%If the analog of Theorem A holds for a single finite set $T$ then the pro-$p$ geometric fundamental group of $\Cal U$ would be finitely generated in this case but it is actually free on an infinite number of generators (Theorem B). 

For example, suppose $\char (K)=0$ and $\char(k)=p>0$. let $\Cal U=\Sp K<T,\frac{1}{T}>$ be the closed annulus of thickness $0$, assume $K$ contains a primitive $p$-th root of unity $\zeta$ and set $\lambda=\zeta-1$. Consider the \'etale $\mu_p$-torsor $f:\Cal V\to \Cal U$ given by the equation $Z^p=1+\lambda ^pT^{-m}$ where $m\ge 1$ is an integer prime-to-$p$.  
Consider $X_K=\Bbb P^1_K$ as the generic fibre of the formal $R$-projective line $X$ obtained by glueing the formal closed discs 
$\Sp R<T>$ and $\Sp R<\frac{1}{T}>$ along the formal annulus $U=\Spf R<T,\frac{1}{T}> $.
Then any Galois extension $\tilde f:Y_K\to X_K$ of $f$ is ramified inside the closed disc $\Sp K<T>$ (which is embedded in $X_K$) above at least $m$ points.
Indeed, let $g:Y\to X$ be the finite generically $\mu_p$-torsor where $Y$ is the normalisation of $X$ in $Y_K$. Then the finite morphism $Y_k\to X_k$ is generically \'etale and 
(generically) defined by an equation $h^p-h=t^{-m}$,
where $h$ (resp. $t$) is the reduction of $H$ defined by $Z=1+\lambda H$ (resp. of $T$) modulo $\pi$, which is a generically \'etale Artin-Schreier cover with conductor $m$ at $0$. 
An easy verification using the Riemman-Hurwitz formula 
(and an argument reducing to the case where the extension $\tilde f$ is \'etale above $\Sp K<\frac{1}{T}>$) 
shows that any Galois extension $\tilde f:Y_K\to X_K$ of $f$ as above is ramified inside the closed disc $\Sp K<T>$ above at least $m$ points. 
As $m$ increases one sees that it is not possible to bound the number of additional branched points in general.
\enddefinition

\definition {Examples 3.9} Suppose $\char(K)=0$ and $\char (k)=p>0$.
Let $U=\Spf R<T>$ (resp. $U=\Spf R<T_1,T_2>/(T_1T_2-1)$) be the standard formal closed unit disc (resp. formal closed annulus of thickness $0$)
embedded in the $R$-projective line $X=\Bbb P^1_{R}$ and $X_K\setminus U_K$ is an open unit disc centred at $\infty$
(resp. embedded in the $R$-formal model of the projective line $\Bbb P^1_K$ consisting
of two standard formal closed unit discs $D_1$ and $D_2$ centred at $0$ and $\infty$; respectively, which are patched with $U$ along their boundaries ($\vert T_1\vert=\vert T_2\vert=1$) 
and $X_K\setminus U_K$ is the disjoint union of two open unit discs).
Let $\Cal U\defeq U_K$ and $S=\{y_1,\cdots,y_r\} \subset \Cal U(K)$ a set of $r\ge 0$ distinct $K$-rational points. 
The results of $\S3$ in this case read as follows. First, the homomorphism $\pi_1(\Cal U\setminus S,\eta)^{\geo}
\to \underset{T} \to {\varprojlim}\  \pi_1(\Bbb P^1_K\setminus (T\cup S),\eta)^{\geo}$, where the projective limit is over all finite sets of points 
$T\subset X_K\setminus \Cal U$, makes $\pi_1(\Cal U\setminus S,\eta)^{\geo}$ 
into a semi-direct factor 
of $\underset{T} \to {\varprojlim}\ \pi_1(\Bbb P^1_K\setminus (T\cup S),\eta)^{\geo}$, the maximal pro-$p$ quotient $\pi_1(\Cal U\setminus S,\eta)^{\geo,p}$
is pro-$p$ free of infinite rank, and the maximal prime-to-$p$ quotient $\pi_1(\Cal U\setminus S,\eta)^{\geo,p'}$
is (pro-)-prime-to-$p$ free of rank $r$ (resp. $r+1$).
\enddefinition

$$\text{References.}$$

\noindent
[Bosch-L\"utkebohmert-Raynaud] Bosch, S., L\"utkebohmert, W., and Raynaud, M. Formal and rigid geometry IV. The reduced fibre theorem, Invent. Math. 119, 
361-398 (1995).

\noindent
[Bourbaki] Bourbaki,N. Commutative Algebra Chapters 1-7, Springer (1989).

\noindent
[Epp] Epp, H. Eliminating wild ramification, Invent. Math. 19, 235-249, (1973).

\noindent
[Garuti] Garuti, M. Prolongements de rev\^etements galoisiens en g\'eom\'etrie 
rigide, Compositio Mathematica, tome 104, n 3 (1996), 305-331.

\noindent
[Grothendieck] Grothendieck, A. Rev\^etements \'etales et groupe fondamental, Lecture 
Notes in Math. 224, Springer, Heidelberg, 1971.

\noindent
[Raynaud] Raynaud, M. Rev\^etements de la droite affine en caract\'eristique $p>0$
et conjecture d'Abhyankar, Invent. math. 116, 425-462 (1994).

\noindent
[Ribes-Zalesskii] Ribes, L., and Zalesskii, P. Profinite groups, Ergebnisse der Mathematik 
und ihrer Grenzgebiete, Folge 3,  Volume 40.

\noindent
[Serre] Serre, J-P. Cohomologie Galoisienne, Lecture Notes in Math., 5, Springer Verlag, 
Berlin, 1994.

\noindent
[Serre1] Serre, J-P. Construction de rev\^etements \'etale de la droite affine en caract\'eristique $p>0$, C. R. Acad. Sci. Paris 311 (1990), 341-346.

\noindent
[Van Der Put] Van Der Put, M. The class group of a one dimensional affinoid space, Annales de l'institut Fourier, tome 30, n$^{o}$ 4(1980), p.155-164.

\bigskip

\noindent
Mohamed Sa\"\i di

\noindent
College of Engineering, Mathematics and Physical Sciences

\noindent
University of Exeter

\noindent
Harrison Building

\noindent
North Park Road

\noindent
EXETER EX4 4QF 

\noindent
United Kingdom

\noindent
M.Saidi\@exeter.ac.uk

\end
\enddocument